\documentclass[11pt,showkywords,showsubjclass]{amsart}

\title[Thick isotopy property and Heegaard splittings]
{Thick isotopy property and the mapping class groups of Heegaard splittings}

\author{Daiki Iguchi}
\address{Higashihiroshima, Japan
}
\email{diguchi00@gmail.com}

\keywords{3-manifold, Heegaard splitting, mapping class group, 
topologically minimal surface}
\subjclass{57M60, 57M50}

\usepackage{amsthm}
\usepackage{mathrsfs}
\usepackage{latexsym}
\usepackage[dvipdfmx]{graphicx}
\usepackage[abs]{overpic}
\usepackage{caption}
\usepackage[all]{xy}
\usepackage[normalem]{ulem}
\usepackage{marginnote}
\usepackage{comment}
\usepackage{cite}
\usepackage{overpic}
\usepackage{bm}
\usepackage{color}

\theoremstyle{plain}
\newtheorem*{theorem*}{Theorem}
\newtheorem*{lemma*} {Lemma}
\newtheorem*{corollary*} {Corollary}
\newtheorem*{proposition*}{Proposition}
\newtheorem*{conjecture*}{Conjecture}
\newtheorem{theorem}{Theorem}[section]
\newtheorem{lemma}[theorem]{Lemma}
\newtheorem{corollary}[theorem]{Corollary}
\newtheorem{proposition}[theorem]{Proposition}

\newtheorem{claim}{Claim}

\theoremstyle{remark}

\newtheorem*{definition}{Definition}
\newtheorem*{claim*}{Claim}

\theoremstyle{definition}

\newtheoremstyle{citing}
  {}
  {}
  {\itshape}
  {}
  {\bfseries}
  {.}
  {.5em}
  {\thmnote{#3}}

\theoremstyle{citing}

\numberwithin{equation}{section}

\newcommand{\MCG}{\mathrm{MCG}}

\newcommand{\Int}{\operatorname{int}}

\newcommand{\He}{\mathcal{H}}

\newcommand{\Diff}{\mathrm{Diff}}
\newcommand{\Sig}{\mathrm{\Sigma}}

\newcommand{\Area}{\mathrm{Area}}
\newcommand{\ol}[1]{\overline{#1}}

\begin{document}

\maketitle
    
\begin{abstract} 
We give a necessary and sufficient condition for the fundamental group of the space of Heegaard splittings of an irreducible $3$-manifold to be finitely generated. 
The condition is exactly the conclusion of the thick isotopy lemma  
proved by Colding, Gabai and Ketover, 
which says that any isotopy of a Heegaard surface 
is achieved by a $1$-parameter family of surfaces with area bounded above 
by a universal constant 
and with some ``thickness property''. 
We also prove that a Heegaard splitting of a hyperbolic or 
spherical $3$-manifold satisfies the condition 
if it is topologically minimal (in the sense of Bachman)  
and its disk complex has finitely generated homotopy group. 
In conclusion, such a Heegaard splitting has finitely generated mapping class group. 
\end{abstract}

\section{Introduction}\label{sec:intro}
Let $M$ be a closed orientable $3$-manifold. 
A {\em Heegaard splitting} is a decomposition of $M$ into two handlebodies
along a closed embedded surface $\Sigma$. 
We will denote such a splitting of $M$ by $(M,\Sigma)$. 
In \cite{JM}, Johnson and McCullough defined the space $\mathcal{H}(M,\Sigma)$ 
of Heegaard splittings equivalent to $(M,\Sigma)$ by $\Diff(M)/\Diff(M,\Sigma)$, 
where $\Diff(M)$ is the space of self-diffeomorphisms of $M$ and 
$\Diff(M,\Sigma)$ is its subspace consisting of maps that send $\Sigma$ to itself. 
For example, computing the $0$-th homotopy group of $\mathcal{H}(M,\Sigma)$ is 
the same as classifying Heegaard splittings up to isotopy. 
Throughout the paper, we will focus only on the case that $M$ is irreducible. 
In \cite{JM}, the $k$-th homotopy group of $\mathcal{H}(M,\Sigma)$ 
was computed for $k \ge 2$. 
On the other hand, $\pi_1(\mathcal{H}(M,\Sigma))$ is closely related to 
the mapping class group of a Heegaard splitting or the Goeritz group, 
and these groups are still mysterious. 
In this paper, we give a necessary and sufficient condition 
for $\pi_1(\mathcal{H}(M,\Sigma))$ to be finitely generated. 

In \cite{CGK18}, Colding, Gabai and Ketover found an effective algorithm to 
construct the complete list of Heegaard splittings of a non-Haken hyperbolic $3$-manifold. 
A key of their argument is the {\em thick isotopy lemma} (\cite[Lemma~2.10]{CGK18}), 
which allows us to turn the computation of 
the $0$-th homotopy group of the space of Heegaard splittings 
into a purely combinatorial problem, 
involving the (crudely) almost normal surface theory. 
The same strategy is also useful in computing $\pi_1(\mathcal{H}(M,\Sigma))$ 
as stated below. 
From now on, we fix a Riemannian metric on $M$. 
Let $\delta>0$. 
A surface $S$ in $M$ is said to be 
{\it $\delta$-compressible} if 
there exists a compressing disk $D$ for $S$ 
such that $\mathrm{diam}\,\partial D \le \delta$. 
Otherwise $S$ is said to be {\it $\delta$-locally incompressible}.   

\begin{definition}
We say $(M,\Sigma)$ satisfies the {\em thick isotopy property} if the following holds. 
There exist $C>0$ and $\delta>0$, depending only on $\Sigma$ and 
the metric of $M$, such that 
any isotopy ${\{\Sigma_{t}\}}_{t \in I}$ with $\Sigma_0=\Sigma_1=\Sigma$ 
can be deformed within its homotopy class (as a loop in $\mathcal{H}(M,\Sigma)$) 
so that afterward for all $t \in I$,  
\begin{itemize}
\item $\mathrm{Area}(\Sigma_t) < C$, and
\item $\Sigma_t$ is $\delta$-locally incompressible.  
\end{itemize}
\end{definition}

\begin{theorem}\label{thm:finitely_generatedness} 
The fundamental group of $\mathcal{H}(M,\Sigma)$ is finitely generated if and only if 
$(M,\Sigma)$ satisfies the thick isotopy property. 
\end{theorem}

Our second aim is to investigate what kind of Heegaard splitting 
satisfies the thick isotopy property.   
Let $S$ be a closed embedded surface in $M$ of genus at least $2$. 
The {\em disk complex} $\Gamma(S)$ of $S$ is defined to be the 
simplicial complex 
whose vertices are the isotopy classes of compressing disks for $S$,  
and whose $i$-simplices are $(i+1)$-tuples of vertices 
that admit disjoint representatives. 
In \cite{Bac10}, Bachman introduced the concept of a topologically minimal surface 
as a generalization of several important classes of surfaces in a $3$-manifold, 
including incompressible surfaces and strongly irreducible surfaces.  

\begin{definition}[Bachman~\cite{Bac10}]
We say $S$ is {\em topologically minimal} if 
$\Gamma(S)=\emptyset$ or $\pi_{d-1}(\Gamma(S)) \neq 1$ 
for some $d \in \mathbb{N}$. 
If $S$ is topologically minimal, 
the {\em topological index} of $S$ is defined to be 
the smallest number $d$ such that $\pi_{d-1}(\Gamma(S)) \neq 1$.  
\end{definition}

\begin{theorem}\label{thm:TIP_minimal_surface}
Let $M$ be a hyperbolic or spherical $3$-manifold but not $S^3$. 
Let $(M,\Sigma) $ be a Heegaard splitting of $M$. 
Suppose that $\Sigma$ is a topologically minimal surface of index $d$. 
Furthermore, suppose that $\pi_{d-1}(\Gamma(\Sigma))$ is finitely generated 
if $d>1$. 
Then $(M,\Sigma)$ satisfies the thick isotopy property.  
\end{theorem}

Here are a few remarks on the theorem. 
First, we note that the index of $\Sigma$ is $1$  
if and only if $\Sigma$ is strongly irreducible, 
and so there are many examples of Heegaard splittings 
satisfying the assumption of the theorem. 
Unfortunately, we do not know if such examples exist when $d>1$.  
(This question can also be seen as the special case of Question~5.10 in \cite{Bac10}.)
However, the examples by Campisi-Rathbun \cite{CR18} 
possibly satisfy the assumption.  
Building on the idea of  Bachman-Johnson \cite{BJ10}, 
they constructed examples of hyperbolic $3$-manifolds 
that contain a Heegaard surface with index $d$ for every $d>0$. 
Indeed, they proved that there exists a retraction from 
the disk complex to a sphere $P \subset \Gamma(\Sigma)$ of appropriate dimension. 
It is likely that we can arrange the construction 
so that such a sphere is in fact a deformation retract. 

We also note that the assumption on $\pi_{d-1}(\Gamma(\Sigma))$ 
in Theorem~\ref{thm:TIP_minimal_surface} 
is used only in the proof of Lemma~\ref{lem:extension},  
and hence it can be replaced with any condition 
that implies Lemma~\ref{lem:extension}. 
In particular, it would be interesting to search for examples of Heegaard splittings 
for which Lemma~\ref{lem:extension} holds without the assumption. 
In Section~\ref{sec:example}, 
we see that this is the case 
for infinitely many examples of Heegaard splittings of $(\mathrm{surface}) \times I$. 
As a consequence, those Heegaard splittings satisfy the thick isotopy property. 

The above theorems have an application 
to the theory of the mapping class group of a Heegaard splitting.  
For a Heegaard splitting $(M,\Sigma)$, 
its {\em mapping class group} $\mathrm{MCG}(M,\Sigma)$  
is defined to be $\pi_0(\Diff(M,\Sigma))$.

\begin{corollary}\label{cor:MCG_minimal surface}
If $M$ and $\Sigma$ are as in Theorem~\ref{thm:TIP_minimal_surface}, then  
$\mathrm{MCG}(M,\Sigma)$ is finitely generated. 
\end{corollary}

In Section~\ref{sec:example}, 
we also establish finite generation of the mapping class groups for 
infinitely many examples of Heegaard splittings of $(\mathrm{surface}) \times I$. 

There have been many efforts to find a finite generating set for 
the mapping class group of a Heegaard splitting. 
Possibly the most interesting is $\MCG(S^3,\Sigma_g)$, 
where $(S^3,\Sigma_g)$ is a standard genus $g$ Heegaard splitting of the $3$-sphere. 
It is known that $\MCG(S^3,\Sigma_g)$  is finitely generated for $g=2$ by 
\cite{Goe33} (see also \cite{Sch04}), and for $g=3$ by \cite{FS18}. 
However, it is not known if the same is true for $g \ge 4$. 
On the other hand,  
a genus $\ge 2$ Heegaard surface in $S^3$ 
is topologically minimal by \cite{App10} or \cite{CT20}.  
(But the disk complex is not of finite type. See Appendix~\ref{sec:appA}.) 
So there might be a good chance to improve our proof 
to remove the assumption that $M \neq S^3$. 
In fact, much of our argument is still valid when $M=S^3$:  
Lemma~\ref{lem:s-map} below is the only place 
where the assumption $M \neq S^3$ is used essentially. 

As another example, 
any genus $2$ weakly reducible Heegaard splitting 
has finitely presented mapping class group by 
\cite{Cho13, CK14, CK15, CK16, CK19,Ak08, Cho08}. 
A finite generating set for a genus $3$ Heegaard splitting of the $3$-torus 
is also known by \cite{Joh11}. 
While little has been known about the mapping class group of a Heegaard splitting 
of genus greater than $3$, 
an advantage of our approach is 
that it is applicable to arbitrarily high genus Heegaard splittings. 

\subsection*{Organization of the paper}
Section~\ref{sec:Normal_surface} is a preliminary towards the proof of 
Theorem~\ref{thm:finitely_generatedness}, including the definition of 
a crudely almost normal surface.  
Theorem~\ref{thm:finitely_generatedness} 
is proved in Section~\ref{sec:proof_fg}. 
Section~\ref{sec:min-max_thm} is a quick introduction to min-max theory. 
In Section~\ref{sec:proof_TIP_minimal_surface}, 
we prove Theorem~\ref{thm:TIP_minimal_surface} and 
Corollary~\ref{cor:MCG_minimal surface}. 
In Section~\ref{sec:example}, 
we prove that  
infinitely many examples of Heegaard splittings of $({\rm surface}) \times I$ 
satisfy the thick isotopy property, 
and as a consequence their mapping class groups are finitely generated. 
In Appendix~\ref{sec:appA}, we see that the disk complex of a Heegaard surface of $S^3$ 
is not homotopy equivalent to a finite simplicial complex.  

\subsection*{Acknowledgements}
The author would like to thank Professor Yuya Koda 
for his valuable comments, advice, and helpful conversations. 

\section{Normal surface theory}\label{sec:Normal_surface}
Throughout the paper, we will use the following notations: 
\begin{itemize}
\item $I:=[0,1]$.  
\item For $r>0$, $B^d_r:=\{x \in \mathbb{R}^d \mid |x| \le r\}$. 
\item If $\mathcal{K}$ is a simplicial complex, we will denote by $\mathcal{K}^i$ 
its $i$-skeleton. 
\end{itemize}

In this section, we recall some definitions and lemmas from \cite{CGK18}.  
Let $M$ be a closed orientable $3$-manifold. 
Let $\mathcal{T}$ be a triangulation of $M$. 

\begin{definition}
A closed embedded surface $S \subset M$ 
is {\it crudely almost normal} (with respect to $\mathcal{T}$) if 
the following are satisfied: 
\begin{enumerate}
\item $S$ is transverse to any simplex of $\mathcal{T}$. 
\item If $\tau$ is a $2$-simplex of $\mathcal{T}$, 
$S \cap \tau$ consists of finitely many arcs (with no circle component). 
\item If $\sigma$ is a $3$-simplex of $\mathcal{T}$,  
$S \cap \sigma$ consists of finitely many disks 
but possibly with one exception: 
there may be exactly one $3$-simplex 
that contains exactly one unknotted annulus component.  
\end{enumerate}
A crudely almost normal surface will be called a {\it crudely normal surface} 
if it has no exceptional annulus component.  
\end{definition}

The {\it weight} of $S$ is defined to be $|S \cap \mathcal{T}^1|$. 
Let $S'$ be another crudely almost normal surface. 
Then, $S$ and $S'$ are said to be {\it normally isotopic} 
if they are isotopic through surfaces transverse to each simplex. 
We say $S'$ is {\it obtained from $S$ by a pinch} if 
it is obtained from $S$ and a $2$-sphere in $M$ 
by connecting them with a tube.  
Such a move or its inverse will be called a {\it pinch}. 
Note that if $M$ is irreducible, 
a pinch can be achieved by an isotopy 
that contracts the $2$-sphere across a $3$-ball in $M$. 

\begin{lemma}[{\cite[Lemma~3.4]{CGK18}}]\label{lem:finite}
There are only finitely many normal isotopy classes 
of crudely almost normal surfaces with weight at most $L$. 
\end{lemma}

\begin{figure}
\begin{overpic}[width=12cm]{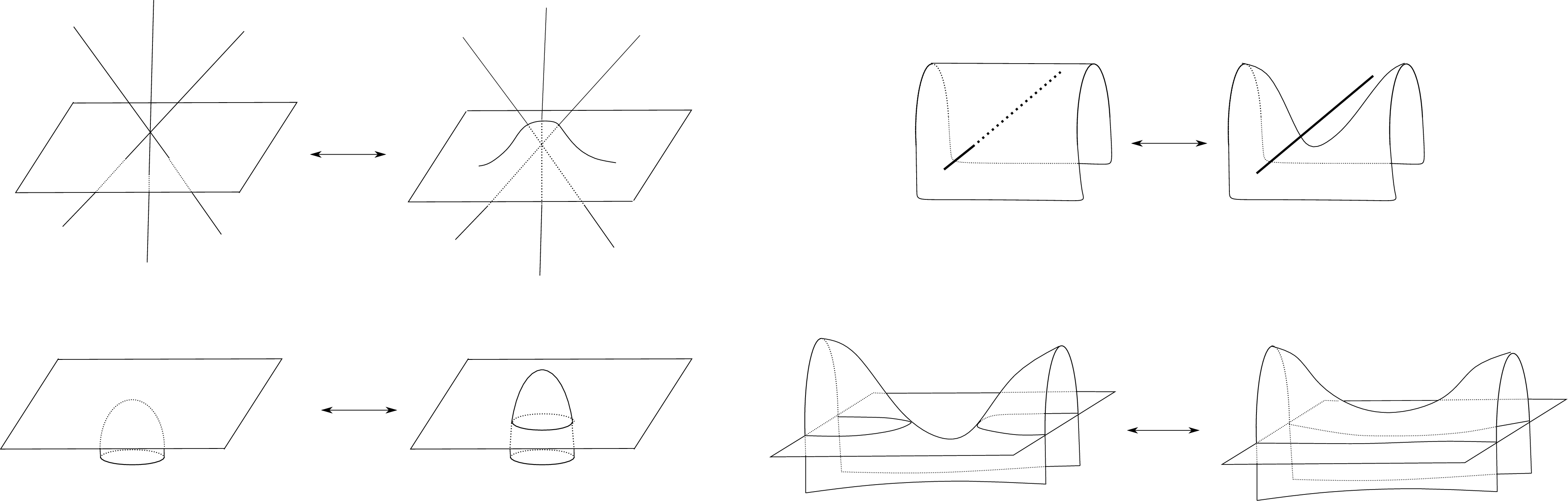}
\put(68,55){(0)}
\put(245,60){(1)}
\put(68,0){(2a)}
\put(241,0){(2b)}
\end{overpic}
\caption{When $t$ moves from $t_i-\epsilon$ to $t_i+\epsilon$, 
$\Sigma_t$ passes either (0) a $0$-simplex, 
(1) a $1$-simplex, 
(2a) a center tangency with a $2$-simplex or (2b) a saddle tangency. 
}\label{fig:tangency}
\end{figure}

A {\em generic $\mathcal{T}$-isotopy}  is an isotopy $\{S_t\}_{t \in I}$ such that 
$S_t$ is transverse to $\mathcal{T}$ for all $t$ but finitely many points 
$0<t_1,\ldots,t_l<1$. 
In addition, $S_{{t_i}+\epsilon}$ and $S_{{t_i}-\epsilon}$ 
differ by one of the moves shown in Figure~\ref{fig:tangency}: 
When $t$ moves from $t_i-\epsilon$ to $t_i+\epsilon$, 
$\Sigma_t$ passes either (0) a $0$-simplex, 
(1) a $1$-simplex, (2a) a center tangency with a $2$-simplex or (2b) a saddle tangency. 
Furthermore, if $|S_t \cap \mathcal{T}^1| \le L$ for $t \in I$, 
$\{S_t\}_{t \in I}$ is called a {\em generic $L$-$\mathcal{T}$-isotopy}. 
Finally, for a Riemannian $3$-manifold $M$, 
an isotopy $\{S_t\}_{t \in I}$ is said to be a 
{\em $C$-isotopy} if $\mathrm{Area}(S_t) \le C$ for $t \in I$. 

\begin{lemma}
[{\cite[Lemma~3.2]{CGK18}}]
\label{lem:weight_bound}
Let $\mathcal{T}$ be a triangulation of the Riemannian $3$-manifold $M$ 
with metric $\rho$, $L>0$ and $\varepsilon>0$. Then there exists 
$K(\mathcal{T},L,\epsilon,\rho)>0$ such that if $S$ is a closed embedded surface with 
$\mathrm{Area}(S)<C$, then $S$ is isotopic to a surface $S'$ such that 
$|S' \cap \mathcal{T}^1|<KC$ and the diameter of the trace of any point of the isotopy is at most 
$\varepsilon$. 

If $F:S \times [0,1] \rightarrow M$ is a $C$-isotopy between surfaces $S_0$ and $S_1$ 
that are transverse to $\mathcal{T}$ of weight at most $L$, then there exists a generic $K(C+1)$-$\mathcal{T}$-isotopy $G$ from $S_0$ to $S_1$ such that, 
for all $x \in S$ and $t \in [0,1]$, $d(G(x,t),F(x,t))<\varepsilon$.
\end{lemma}

\section{The proof of Theorem~\ref{thm:finitely_generatedness}}\label{sec:proof_fg}
In this section, we prove Theorem~\ref{thm:finitely_generatedness}.  
It is not hard to see the necessity, i.e. forward implication, in the theorem.  
First, recall that a path in $\He(M,\Sigma)$ 
can be identified with an isotopy of a Heegaard surface. 
More precisely, 
as $\Diff(M) \to \Diff(M)/\Diff(M,\Sigma)=\He(M,\Sigma)$ is a fibration \cite{JM}, 
any path $\alpha:I \to \He(M,\Sigma)$ lifts to $\tilde{\alpha} :I \to \Diff(M)$ 
and we can define an isotopy of a Heegaard surface 
by $\Sigma_t:=\tilde{\alpha}(t)(\Sigma)$. 
Conversely, if an isotopy of a Heegaard surface is given, 
it defines a path in $\He(M,\Sigma)$ via the isotopy extension theorem. 
Now if $\pi_1(\He(M,\Sigma))$ is finitely generated, 
we can find a finite collection of isotopies of $\Sigma$ 
such that any isotopy representing an element of $\pi_1(\He(M,\Sigma))$ 
can be expressed as 
the product of isotopies in the collection. 
Thus, $(M,\Sigma)$ satisfies the thick isotopy property. 

In the following, 
we prove the sufficiency of the theorem. 
Let $C>0$ and $\delta>0$ be the constants given 
in the definition of the thick isotopy property. 
Fix a triangulation $\mathcal{T}$ of $M$ such that 
\begin{itemize}
\item $\Sigma$ is crudely normal with respect to $\mathcal{T}$, and 
\item any simplex of $\mathcal{T}$ has the diameter at most $\delta$.  
\end{itemize}
Let $\{\Sigma_t\}_{t \in I}$ be any isotopy with $\Sigma_0=\Sigma_1=\Sigma$. 
By the argument in \cite{CGK18}, we can convert $\Sigma_t$ to a crudely almost normal surface with respect to $\mathcal{T}$. 
Here is a very rough sketch of the argument. 
By assumption, $\Area(\Sigma_t) \le C$ and $\Sigma_t$ is $\delta$-locally incompressible for $t \in I$. By Lemma~\ref{lem:weight_bound}, $\Sigma_t$ is transverse to every simplex of $\mathcal{T}$ for all but finitely many points and $\Sigma_t$ has weight at most $L:=K(C+1)$. Using the $\delta$-locally incompressibility condition, for every $3$-simplex $\sigma$ of $\mathcal{T}$, we can pinch off and remove non-disk components of $\Sigma_t \cap \sigma$. (We note that a subtle situation may occur when $\Sigma_t$ passes through a tangency of type (2b). Around such a tangency, we may be forced to allow an unknotted annulus component. See \cite[Proof of Lemma~3.6]{CGK18} for more details.)  
In summary, we have 

\begin{claim}\label{clm:normalization}
$\{\Sigma_t\}_{t \in I}$ can be deformed 
within its homotopy class (as the loop in the space $\mathcal{H}(M,\Sigma)$) 
so that afterward for all $t$ 
but finitely many points in $I$, $\Sigma_t$ is 
a crudely almost normal surface with weight bounded above 
by a universal constant $L>0$.  
\end{claim}

\begin{proof}
This follows from \cite[Lemma~3.6]{CGK18}. 
\end{proof}

Consider the graph $\mathcal{G}$ such that 
each vertex of $\mathcal{G}$ corresponds to a normal isotopy class 
of crudely almost normal surfaces w.r.t. $\mathcal{T}$ with weight at most $L$, and 
each edge corresponds to one of the moves~(0) -- (2b) 
shown in Figure~\ref{fig:tangency} or a pinch. 
By Lemma~\ref{lem:finite}, $\mathcal{G}$ is a finite graph. 
In particular, $\pi_1(\mathcal{G})$ is finitely generated. 
By Claim~\ref{clm:normalization}, the natural homomorphism 
$\pi_1(\mathcal{G}) \to \pi_1(\mathcal{H}(M,\Sigma))$ is a surjection. 
Thus, we conclude that $\pi_1(\mathcal{H}(M,\Sigma))$ is finitely generated. \qed

\section{The min-max theorem}\label{sec:min-max_thm}
This section is a quick introduction to the min-max theory of Simon-Smith \cite{Sm82}, 
which will be used in the next section. 
One can consult e.g. \cite{CD03, CGK22} for more details on this subject. 

Let $M$ be a closed, orientable, Riemannian $3$-manifold. 
We will denote by $\mathscr{H}^2(\cdot)$ the $2$-dimensional Hausdorff measure on $M$. 

\begin{definition}
Let $X^k$ be a manifold. 
A family $\{\Sig_t\}_{t \in X}$ of closed subsets of $M$ is called a 
{\it (genus $g$) sweep-out} if it satisfies the following conditions: 
\begin{itemize}
\item $\Sig_t$ converges to $\Sig_{t_0}$ in the Hausdorff topology when $t \to t_0$. 
\item $\mathscr{H}^2(\Sig_t) \rightarrow \mathscr{H}^2(\Sig_{t_0})$ 
	when $t \rightarrow t_0$. 
\item $\Sigma_t$ is a closed genus $g$ surface in $M$ if $t \in \mathrm{int}\,X$. 
On the other hand, if $t \in \partial X$, $\Sigma_t$ is a closed surface 
of $\mbox{genus} \le g$ plus finitely many arcs. 
\item $\Sig_t$ varies smoothly for $t \in \mathrm{int}\,X$. 
\end{itemize}
\end{definition}

For later use, we restrict ourselves to the case that $X=I \times B^d$.  
Consider the subspace $\mathscr{I}$ of 
$C^\infty(M \times (I \times B^d),M)$ consisting of those maps $\psi$ 
such that  
\begin{enumerate}
\item[($\mathrm{i}$)] $\psi(\cdot,t)$ is a diffeomorphism of $M$ 
for $t \in I \times B^d$, and 
\item[($\mathrm{ii}$)] $\psi(\cdot,t)=\mathrm{id}_M$ for $t \in \partial I \times B^d$.  
\end{enumerate} 
Let $\mathscr{I}_0 \subset \mathscr{I}$ 
be the component containing the map $\psi_0$ given by 
$\psi_0(x,t):=x$. 
Given a sweep-out $\{\Sig_t\}_{t \in I \times B^d}$, 
define the collection $\Pi_{\{\Sig_t\}}$ of sweep-outs by 
\[\Pi_{\{\Sig_t\}}:=
\left\{ \psi(\Sigma_t,t)\}_{t \in I \times B^d} 
\mid \psi \in \mathscr{I}_0
\right\}.
\] 
The {\it width} of $\Pi_{\{\Sig_t\}}$ is defined by 
\[W(\Pi_{\{\Sig_t\}},M):=
\inf_{\{\Lambda_t\} \in \Pi_{\{\Sig_t\}}}\sup_{t \in I \times B^d}\mathscr{H}^2(\Lambda_t).\]

A sequence $\{\Sig_t^i\}_{t \in X}$ ($i \in \mathbb{N}$) of sweep-outs 
in $\Pi_{\{\Sig_t\}}$ is 
a {\it minimizing sequence} if 
$W(\Pi_{\{\Sig_t\}},M)=\lim_{i \to \infty}\sup_{t \in I \times B^d}\mathscr{H}^2(\Sig_t^i)$. 
Furthermore, a sequence $\{\Sig_{t_i}^i\}_{i \in \mathbb{N}}$ is 
a {\it min-max sequence} if 
$W(\Pi_{\{\Sig_t\}},M)=\lim_{i \to \infty}\mathscr{H}^2(\Sig_{t_i}^i)$. 

Simon-Smith's min-max theorem is the following theorem. 
(The following statement can be found in \cite{CGK18} with minor modification, 
see \cite{CD03,DP10,Ke19} for the proof 
and also \cite[Appendix]{CGK18} for the multi-parameter case.)  

\begin{theorem}
[cf. {\cite[Theorem~2.1]{CGK18}}]
\label{thm:min-max} 
Given a sweep-out $\{\Sig_t\}_{t \in I \times B^d}$ of genus-$g$ surfaces, if 
\begin{equation}
W(\Pi_{\{\Sig_t\}},M)>\sup_{t \in \partial I \times B^d} \mathscr{H}^2(\Sig_t),
\end{equation}
then there exists a min-max sequence $\Sig_i:=\Sig_{t_i}^i$ such that 
\begin{equation}
\Sig_i \rightarrow \sum_{i=1}^k n_i\Gamma_i \mbox{\hspace{0.5em}as varifolds}, 
\end{equation}
where $\Gamma_i$ are smooth closed embedded minimal surfaces and $n_i$ are 
positive integers. 
Moreover, after performing finitely many compressions on $\Sig_i$ 
and discarding some components, 
each connected component is isotopic to one of the $\Gamma_i$ or to a 
double cover of one of the $\Gamma_i$. 
We have the following genus bounds with multiplicity: 
\begin{equation}\label{eq:genus}
\sum_{i \in \mathscr{O}}n_i g(\Gamma_i) + 
\frac{1}{2}\sum_{i \in \mathscr{N}}n_i(g(\Gamma_i)-1)
\le g, 
\end{equation}
where $\mathscr{O}$ denotes the subcollection of $\Gamma_i$ that is orientable 
and $\mathscr{N}$ denotes those $\Gamma_i$ 
that are nonorientable, 
and where $g(\Gamma_i)$ denotes the genus of $\Gamma_i$ 
if it is orientable, and 
the number of cross-caps that one attaches to a sphere to obtain a homeomorphic 
surface if $\Gamma_i$ is nonorientable.   
\end{theorem}

\begin{lemma}
\label{lem:area_bound}
Let $M$ be a hyperbolic or spherical $3$-manifold. 
If $\{\Sigma_t\}_{t \in I \times B^d}$ is a genus $g$ sweep-out satisfying 
$W(\Pi_{\{\Sig_t\}},M)>\sup_{t \in \partial I \times B} \mathscr{H}^2(\Sig_t)$, 
then $W(\Pi_{\{\Sig_t\}},M) \le 8\pi(g+1)$. 
\end{lemma}

\begin{proof}
If $M$ is hyperbolic, 
$W(\Pi_{\{\Sig_t\}},M) \le 4\pi(g-1)$ by \cite[Lemma~9.4]{CGK22}. 
So, we prove the lemma when $M$ is spherical. 
Theorem~\ref{thm:min-max} shows that 
$W(\Pi_{\{\Sig_t\}},M)=\sum_{i=1}^k n_i \mathrm{Area}(\Gamma_i)$ for 
some embedded minimal surfaces $\Gamma_i$ ($1 \le i \le k$). 
By Frankel's theorem \cite{Fra66}, the min-max limit is in fact connected 
and we can express the width as $W(\Pi_{\{\Sig_t\}},M)= n \mathrm{Area}(\Gamma)$. 
If $\Gamma$ is non-orientable, its double cover is 
stable by \cite[Theorem~7.2]{WZ23}. 
But this is impossible because $S^3$ with the standard metric and thus its quotient cannot contain a stable minimal surface.  
So $\Gamma$ must be orientable and again by \cite[Theorem~7.2]{WZ23}, 
the multiplicity $n$ must be one. 
This together with Choi-Schoen's area bound \cite{CS85}
for a minimal surface in a spherical $3$-manifold implies 
\[
W(\Pi_{\{\Sig_t\}},M) =
\mathrm{Area}(\Gamma) \le 
8\pi \left(\frac{2}{|\pi_1(M)|}-\frac{\chi(\Gamma)}{2}\right) \le 
8\pi(g+1). 
\]
\end{proof}

\section{Proof of Theorem~\ref{thm:TIP_minimal_surface}}\label{sec:proof_TIP_minimal_surface}

In this section, we prove Theorem~\ref{thm:TIP_minimal_surface}. 
Suppose that $M$ is a hyperbolic or spherical $3$-manifold but not $S^3$, 
and $\Sigma$ is a genus $g$ Heegaard surface of topological index $d>0$. 
Furthermore, if $d>1$, we assume that $\pi_{d-1}(\Gamma(\Sigma))$ is finitely generated. 

\subsection*{Step~1: The definitions of $C$ and $\delta$} 
Let $\{\varphi_k:S^{d -1}\to \Gamma(\Sigma) \mid k=1,\ldots, n\}$ 
be a collection of maps that represents a finite generating set 
of $\pi_{d-1}(\Gamma(\Sigma))$. 
Put $B=B^d_1$. 
For each $\varphi_k$, 
we consider the sweep-out 
$\{\Sigma^k_s\}$ parametrized by $B$ and given below.   
(The same construction can be found in \cite{Bac10}.) 

Take a triangulation $\mathcal{K}$ of $S^{d-1}$ for which $\varphi_k$ is simplicial. 
For each vertex $s$ of $\mathcal{K}$, 
choose a representative $D_s$ of $\varphi_k(s) \in \Gamma(\Sigma)$ 
so that if $s$ and $s'$ are in the same simplex of $\mathcal{K}$, 
$D_s$ and $D_{s'}$ are disjoint. 
Let $\widetilde{\mathcal{K}}$ be the cone over $\mathcal{K}$. 
So, $\widetilde{\mathcal{K}}$ is a triangulation of $B$.  
Roughly, $\{\Sigma^k_s\}_{s \in B}$ is defined as follows. 
For the center $0$ of $B$, set $\Sigma^k_0=\Sigma$. 
For each vertex $s$ of $\mathcal{K}$, 
let $\Sigma^k_s$ be the result of compressing $\Sigma$ across $D_s$. 
Furthermore, for a barycenter $s$ of an $l$-simplex of $\mathcal{K}$ 
that is spanned by vertices $s_0,\ldots,s_l \in \mathcal{K}^0$, 
define $\Sigma^k_s$ to be the result of compressing $\Sigma$ simultaneously 
across $D_{s_i}$.  
For general $s \in B$,   
we define $\Sigma^k_s$ as a linear interpolation of the above construction.  

Here is a more formal definition. 
Let $\tilde{\sigma}$ be a $d$-simplex of $\tilde{\mathcal{K}}$ and 
let $0=s_0,s_1,\ldots,s_d$ be the vertices of $\tilde{\sigma}$. 
Denote by $\sigma$ the $(d-1)$-simplex spanned by the vertices $s_1,\ldots,s_d$. 
For $1 \le i \le d$, 
fix a cylinder $N(D_{s_i}) \cong B^2_1 \times [-1,1]$ in $M$ 
such that $B^2_1  \times 0=D_{s_i}$ 
and $N(D_{s_i}) \cap \Sigma=\partial B^2_1 \times [-1,1]$. 
Furthermore, $N(D_{s_i})$ can be chosen to be mutually disjoint. 
Let $\rho_\tau:B^2_1 \times [-1,1] \to B^2_1 \times [-1,1]$ be the map given by 
\[\rho_\tau(x,u)=((1-\tau b(u)) x,u)\]
for $(x,u) \in B^2_1 \times [-1,1]$ and $\tau \in I$. 
Here, $b(u)$ is a bump function. 
Thus, $\rho_\tau$ shrinks the annulus $\partial B^2_1 \times [-1+\epsilon,1-\epsilon]$ 
to an arc $0 \times  [-1+\epsilon,1-\epsilon]$. 
Fix an identification $\iota_{s_i}:B^2_1 \times [-1,1] \to N(D_{s_i})$ for $1 \le i \le d$,  
and define the homotopy $\eta^{s_i}_\tau:M \to M$ to be 
\[\eta^{s_i}_\tau:=\iota_{s_i} \circ \rho_\tau \circ \iota^{-1}_{s_i}\] 
on $N(D_{s_i})$ and 
to be the identity on the complement of a small neighborhood of $N(D_{s_i})$. 
Regard $\tilde{\sigma}$ as a $d$-cube in $\mathbb{R}^d$ 
so that $0=s_0$, the $i$th unit vector represents $s_i$,  
and the other corners correspond to the barycenters 
of the simplices contained in $\sigma$. 
Then, every point $s \in \tilde{\sigma}$ can be written as 
$s=\tau_1s_1+\cdots+\tau_ds_d$ for  
$\tau_1,\ldots,\tau_d \in I$. 
Now we define $\Sigma^k_s$ as the image of $\Sigma$ by the composition 
\[\eta^{s_1}_{\tau_1} \circ \cdots \circ \eta^{s_d}_{\tau_d}.\]
($\eta^{s_i}_{\tau_i}$'s are commutative with each other 
since their supports are mutually disjoint.)  
Since the above construction is consistent 
on every intersection  between adjacent simplices, 
we obtain the sweep-out $\{\Sigma^k_s\}_{s \in B}$. 

For simplicity, given a surface $T$ in $M$, we define 
\[\gamma(T):=
\min \{ \mathrm{diam}\, \partial D \mid \mbox{$D$ is a compressing disk for $T$}\}.\] 
We write $\mathrm{inj}(M)$ for the injectivity radius of $M$. 
Now define 
\[C:=\max \left\{\max_{s \in B}\mathscr{H}^2(\Sigma^1_s), \ldots, 
\max_{s \in B}\mathscr{H}^2(\Sigma^n_s), 
8\pi(g+1) \right\}+1,\] 
and 
\[\delta:=
\min \{3^{-(1+2+\cdots+d+(d+1))}\cdot \mathrm{inj}(M), \gamma(\Sigma)/2\}.\]

\subsection*{Step~2: Lemmas}
Let $k \in \{1,\ldots,n\}$.  
Note that by the isotopy extension theorem, 
we can fix a diffeomorphism $f_s:\Sigma \to \Sigma^k_s$  
for $s \in \mathrm{int}\,B$ simultaneously. 
In the next lemma, which follows from the definition of $\{\Sigma^k_s\}_{s \in B}$, 
we identify $\Gamma(\Sigma)$ with $\Gamma(\Sigma^k_s)$ 
through these diffeomorphisms.  
\begin{lemma}\label{lem:h-nontrivial}
Let $S=\partial B^d_{1-\epsilon}$. 
If $\epsilon$ is small enough, $\Sigma^k_s$ is $\delta$-compressible for $s \in S$.  
Furthermore, if $\mathcal{S}$ is a triangulation of $S$ 
such that the diameter of each simplex is small enough, and 
if we define $\psi:\mathcal{S}^0 \to \Gamma(\Sigma)^0$ 
by sending $s \in \mathcal{S}^0$ to one of $\delta$-compressing disks for $\Sigma_s$ 
(and applying $f_s^{-1}$),  
then $\psi$ determines the simplicial map $\mathcal{S} \to \Gamma(\Sigma)$ 
homotopic to $\varphi_k$. 
\end{lemma}

Define $U_0 \subset B$ to be the set of points $s$ 
such that $\Sigma^k_s$ is $\delta$-compressible. 

\begin{lemma}\label{lem:star}
$B \setminus U_0$ is a star-shaped region. 
\end{lemma}

\begin{proof}
The proof is by contradiction.  
Suppose that 
there exists a line segment $\ell \subset B$ 
connecting $0$ with $s \in B \setminus U_0$ 
such that $\ell$ contains a point $u \in U_0$. 
We can find a compressing disk $D$ for $\Sigma^k_u$ with 
$\mathrm{diam}\,\partial D \le \delta$. 
Let $Q$ be a $3$-ball of diameter $\le \delta$ that contains $\partial D$. 
After perturbing $Q$, $\partial Q$ intersects $\Sigma^k_0=\Sigma$, 
$\Sigma^k_s$ and $\Sigma^k_u$ transversely. 
If $Q \cap \Sigma_s$ contains circles that are essential in $\Sigma^k_s$, 
one of such circles bounds a compressing disk for $\Sigma^k_s$. 
As $\Sigma^k_s$ is $\delta$-locally incompressible, this case cannot occur. 
Thus, all the circles in $Q \cap \Sigma^k_s$ are inessential in $\Sigma^k_s$. 
Similarly, all the circles in $Q \cap \Sigma$ are inessential in $\Sigma$. 
By the innermost disk argument, we can isotope $Q$ so that $Q$ is contained 
in the region between $\Sigma$ and $\Sigma^k_s$, 
which is diffeomorphic to $\Sigma \times I$. 
Note that $Q$ still intersects $\Sigma^k_u$ 
so that one of the circles in $Q \cap \Sigma^k_u$ is  essential in $\Sigma^k_u$. 
Thus, $\Sigma^k_u$ is compressible in the product region. 
This is a contradiction because $\Sigma^k_u$ is isotopic to a level surface 
in the product region.  
\end{proof}

\subsection*{Step~3: Extending the sweep-out}
Fix a nontrivial map $\varphi:S^{d-1} \to \Gamma(\Sigma)$ once.  
We may assume that $\varphi=\varphi_1$ and 
set $\{\Sigma_{0s}\}_{s \in B}:=\{\Sigma^1_s\}_{s \in B}$. 

Let $h_t:\Sigma \to M$ be any isotopy with $h_0(\Sigma)=h_1(\Sigma)=\Sigma$. 
To prove Theorem~\ref{thm:TIP_minimal_surface}, 
we must show that there exists an isotopy $h'_t$ 
equivalent to $h_t$ such that for $t \in I$ 
\begin{itemize}
\item $\mathrm{Area}(h'_t(\Sigma)) < C$, and 
\item $h'_t(\Sigma)$ is $\delta$-locally incompressible. 
\end{itemize}
By the isotopy extension theorem, $h_t:\Sigma \to M$ extends 
to $\tilde{h}_t:M \to M$.  
Now define $\Sigma_{ts}$ for $t \in I$ and $s \in B$ by 
\[\Sigma_{ts}:=\tilde{h}_t(\Sigma_{0s}).\] 

\begin{lemma}\label{lem:extension}
$\{\Sigma_{ts}\}_{(t,s) \in I \times B}$ extends 
to $\{\Sigma_{ts}\}_{(t,s) \in [0,2] \times B}$ such that 
\begin{enumerate}
\item $\Sigma_{t0}=\Sigma$ for $t \in [1,2]$,  
\item $\mathrm{Area}(\Sigma_{2s}) < C$ for $s \in B$, 
\item $B \setminus U_2$ is a star-shaped region, 
where $U_2$ is the set of those points $s$ 
such that $\Sigma_{2s}$ is $\delta$-compressible. 
\end{enumerate}
\end{lemma}

\begin{proof}
First, note that there is a natural action of $\Diff(M,\Sigma)$ on $\Gamma(\Sigma)$, 
which induces the action on $[S^{d-1},\Gamma(\Sigma)]$. 
By construction, $\{\Sigma_{1s}\}_{s \in B}$ is 
the image of $\{\Sigma_{0s}\}_{s \in B}$ by $\tilde{h}_1$. 
In other words, $\{\Sigma_{1s}\}_{s \in B}$ can be recovered 
from $\tilde{h}_1 \cdot \varphi$ as follows. 
We repeat the same construction as in Step~1. 
Take a triangulation $\mathcal{K}$ of $S^{d-1}$ 
such that $\varphi:S^{d-1} \to \Gamma(\Sigma)$ is simplicial. 
For each vertex $s$ of $\mathcal{K}$, 
choose a compressing disk $D_s$ for $\Sigma$ 
that represents $\varphi(s)$. 
Then, $\tilde{h}_1(D_s)$ is a compressing disk 
that represents $\tilde{h}_1 \cdot \varphi(s)$. 
Fix an identification 
$\kappa_s:B^2_1 \times [-1,1] \to \tilde{h}_1(N(D_s))$ for each $s \in \mathcal{K}^0$, 
and define the homotopy $\theta^s_\tau:M \to M$ to be 
\[\theta^s_\tau:=\kappa_s \circ \rho_\tau \circ \kappa_s^{-1}\] 
on $\tilde{h}_1(N(D_s))$ and 
to be the identity on the complement of a small neighborhood of $\tilde{h}_1(N(D_s))$. 
On each simplex $\tilde{\sigma}$ of $\widetilde{\mathcal{K}}$ 
with vertices $0=s_0,s_1,\ldots,s_d$, 
define $\Sigma_{1s}$ by the image of the homotopy 
\[\theta^{s_1}_{\tau_1} \circ \cdots \circ \theta^{s_d}_{\tau_d}:M \to M,\]
where $s=\tau_1s_1+\cdots+\tau_ds_d$ and 
$\tau_1,\ldots,\tau_d \in I$.
Note that the above construction contains an ambiguity 
regarding the choice of an identification 
$B^2_1 \times [-1,1] \to \tilde{h}_1(N(D_s))$:
we have to see that $\kappa_s$ coincide with the image of the ``standard'' one, 
that is $\tilde{h}_1 \circ \iota_s$. 
But these two maps are isotopic 
by the uniqueness of tubular neighborhoods, 
and hence the corresponding sweep-outs can be interpolated via this isotopy. 
In this way, we recover the sweep-out $\{\Sigma_{1s}\}_{s \in B}$. 

By assumption, there is a homotopy 
$\Phi:[1,2] \times S^{d-1} \to \Gamma(\Sigma)$ such that 
$\Phi_1=\tilde{h}_1 \cdot \varphi$ and $\Phi_2$ is a product of $\varphi_k$'s.   
By the relative simplicial approximation theorem (see e.g. \cite{Zee64}), 
we can extend $\{1\} \times \mathcal{K}$ 
to a triangulation $\mathcal{L}$ of $[1,2] \times S^{d-1}$ ($=[1,2] \times \partial B$) 
such that $\Phi$ is simplicial with respect to $\mathcal{L}$. 
In the following, 
we will construct a sweep-out $\{\Sigma_{ts}\}_{(t,s) \in [1,2] \times B}$ 
that ``shadows'' the homotopy $\Phi$. 

The construction is similar to that in Step~1. 
Extend the triangulation $\mathcal{L}$ over $\partial ([1,2] \times B)$ 
by adding the two points $(1,0)$ and $(2,0)$ as vertices. 
Regarding $[1,2] \times B$ as a cone over $\partial ([1,2] \times B)$, 
we obtain a triangulation $\widetilde{\mathcal{L}}$ of $[1,2] \times B$. 
For each vertex $v=(t,s)$ on $\mathcal{L}$, 
choose a representative $D_v$ of $\Phi(v)$ 
so that if $v$ and $v'$ are in the same simplex, 
$D_v$ and $D_{v'}$ are disjoint. 
For $v=(1,0)$ or $(2,0)$, let $D_v=\emptyset$.

As before, fix an identification 
$\kappa_v:B^2_1 \times [-1,1] \to N(D_v)$ 
for each vertex $v$ of $\mathcal{L}$, and 
define $\theta^v_\tau: M \to M$ to be $\kappa_v \circ \rho_\tau \circ \kappa^{-1}_v$ 
on $N(D_v)$ and to be the identity 
outside a small neighborhood of $N(D_v)$. 
It suffices to construct the sweep-out on every simplex $\tilde{\sigma}$ 
of $\widetilde{\mathcal{L}}$. 
Let $v_0,\ldots,v_{d+1}$ be the vertices of $\tilde{\sigma}$ and 
regard $\tilde{\sigma}$ as a $(d+1)$-cube in Euclidean space 
with $v_0=0$. 
For $x=\tau_1v_1+\cdots+\tau_{d+1}v_{d+1} \in \tilde{\sigma}$, 
$\tau_1,\ldots,\tau_{d+1} \in I$, 
define $\Sigma_x$ as the image of $\Sigma$ by the composition 
\[\theta^{v_1}_{\tau_1} \circ \cdots \circ \theta^{v_{d+1}}_{\tau_{d+1}}.\]

Note that for $v \in \{2\} \times S^{d-1}$, 
the disk $D_v$ representing $v$ may not coincide 
with the standard one which has been chosen in Step~1.    
Of course, $D_v$ is isotopic to the standard choice, and  
in addition if $v_1,\ldots,v_l$ be vertices adjacent to $v$ 
and all the $D_{v_i}$ have already been in the standard position, 
then we can isotope $D_v$ to the standard position in the complement of $D_{v_i}$. 
Thus, after these isotopies, we may assume that all the $D_v$ coincide 
with the standard one.  
By the uniqueness of tubular neighborhoods, 
we can make $N(D_v)$ and $\theta^v_\tau$ standard as well. 

This sweep-out $\{\Sigma_{ts}\}_{(t,s) \in [0,2] \times B}$ satisfies the desired property. 
Indeed, if $\ell \subset B$ is a radial segment in $B$, 
$\{\Sigma_{2s}\}_{s \in \ell}$ appears in some $\{\Sigma^k_s\}_{s \in B}$ as a subfamily. 
Thus (2) and (3) hold. 
(1) is obvious from the construction. 
\end{proof}

By Theorem \ref{thm:min-max},  we have the following lemma. 
\begin{lemma}\label{lem:area-bound2}
$\{\Sigma_{ts}\}_{(t,s) \in [0,2] \times B}$ can be modified so that 
afterward $\mathscr{H}^2(\Sigma_{ts}) < C$ for $t \in [0,2]$ and $s \in B$. 
\end{lemma}

\subsection*{Step~4: Lifting a submanifold of $[0,2] \times B$ to $\Gamma(\Sigma)$} 
We fix some notation and terminology.  
Let us fix diffeomorphisms $f_{ts}:\Sigma_{00} \to \Sigma_{ts}$ for 
$(t,s) \in I \times \mathrm{int}\,B$ simultaneously 
via the isotopy extension theorem.  
We can identify $\Gamma(\Sigma_{ts})$ with $\Gamma(\Sigma_{00})$ 
($=\Gamma(\Sigma)$) 
through this identification. 
We say $\Sigma_{ts}$ and $\Sigma_{t's'}$ are {\em $\epsilon$-close} if 
$d_M(f_{ts}(x),f_{t's'}(x))<\epsilon$ for any $x \in \Sigma_{00}$. 
Finally, define $U \subset [0,2] \times B$ to be the set of all points $(t,s)$ 
such that $\Sigma_{ts}$ is $\delta$-compressible. 
 
\begin{lemma}\label{lem:s-map}
Let $Y$ be a $d$-manifold embedded in $U$. 
Let $\epsilon < \delta$. 
Suppose that $\mathcal{Y}$ is a triangulation of $Y$ such that 
if $y,y'$ are in the same simplex of $\mathcal{Y}$, 
then $\Sigma_y$ is $\epsilon$-close to $\Sigma_{y'}$. 
Let $\mathcal{Y}'$ be the barycentric subdivision of $\mathcal{Y}$. 
If we define the map $\psi:\mathcal{Y}^0 \to \Gamma(\Sigma)^0$ 
by sending $y \in \mathcal{Y}^0$ to one of $\delta$-compressing disks 
for $\Sigma_y$ and applying $f_y^{-1}$, 
it extends to a simplicial map $\bar{\psi}:\mathcal{Y}' \to \Gamma(\Sigma)$. 
\end{lemma}

\begin{proof}
To extend $\psi$ to a simplicial map 
$\bar{\psi}: \mathcal{Y}' \to \Gamma(\Sigma_{00})$,  
it suffices to find a collection $\{D_y \mid y \in \mathcal{Y}'^0\}$ of disks 
with the following property.  
\begin{itemize}
\item For $y \in \mathcal{Y}'^0$, 
$D_y$ is a compressing disk for $\Sigma_y$. 
\item If $y$ and $y'$ are in the same simplex of $\mathcal{Y}'$, 
then $f_y^{-1}(D_y)$ and $f_{y'}^{-1}(D_{y'})$ are disjoint 
(i.e. they span a $1$-simplex in $\Gamma(\Sigma_{00})$). 
\end{itemize}
Indeed, if such a collection of disks exists, 
we can define $\bar{\psi}$ by $\bar{\psi}(y)=f_y^{-1}(D_y)$. 
The proof is by induction: we will show 
\begin{claim*}
If $D_{y'}$ has already been defined 
for $y'$ the barycenter of any $(i-1)$-simplex of $\mathcal{Y}$ and 
$\mathrm{diam}\,\partial D_{y'} < 3^{1+ \cdots+i}\cdot \delta$ holds,  
then we can find $D_y$ 
for $y$ the barycenter of any $i$-simplex $\mathcal{Y}$ such that 
$\mathrm{diam}\, \partial D_y < 3^{1+ \cdots+i+(i+1)}\cdot \delta$. 
\end{claim*}

Let $\sigma$ be an $i$-simplex of $\mathcal{Y}$ and let $y$ be the barycenter 
of $\sigma$. 
Let $y_1,\ldots,y_{2^{i+1}-1}=y$ be the vertices of $\mathcal{Y}'$ 
that are contained in $\sigma$. 
By induction, for $1 \le j \le 2^{i+1}-2$, $D_{y_j}$ has already been defined  
and the diameter of $\partial D_{y_j}$ is less than $3^{1+ \cdots+i}\cdot \delta$. 
As $\Sigma_{y_j}$ and $\Sigma_y$ are $\epsilon$-close, 
the image of $\partial D_{y_j}$ on $\Sigma_y$ 
has diameter less than $3^{1+ \cdots+i}\cdot \delta+2\epsilon$. 
In what follows, we work on a single surface, say $\Sigma_y$,  
rather than multiple surfaces. 
We will not distinguish between $D_{y_j}$ and its image on $\Sigma_y$ from their notation. 

After relabeling $y_j$'s if necessary, 
we can assume 
that there exists a number $k$ ($1 \le k \le 2^{i+1}-2$) satisfying the following: 
there exists a metric ball $Q$ with 
$\mathrm{diam}\,Q < k (3^{1+\cdots+i} \cdot \delta+2\epsilon)$ 
such that $\bigcup_{j=1}^k \partial D_{y_j}$ is contained in $Q$ while 
$\bigcup_{j=k+1}^{2^{i+1}-2} \partial D_{y_j}$ is in the complement of $Q$. 
Note that $Q$ is a genuine $3$-ball because 

\begin{align*}
\mathrm{diam}\,Q
&< k (3^{1+\cdots+i} \cdot \delta +2\epsilon)\\
&\le (2^{i+1}-2) \cdot (3^{1+\cdots+i} \cdot \delta +2\delta) \\
&< 3^{1+\cdots+(i+1)}\cdot \delta \\
&\le \mathrm{inj}(M). 
\end{align*}

After perturbing $Q$, we assume 
that $\partial Q$ intersects $\Sigma_y$ transversely. 
We can find a circle in 
$\partial Q \cap \Sigma_y$ that is essential in $\Sigma_y$. 
Indeed, if all the circles in $\partial Q \cap \Sigma_y$ were inessential, 
by the innermost disk argument, 
$\Sigma_y$ could be isotoped so that afterward $\Sigma_y \subset Q$. 
This is impossible 
because $\Sigma_y$ is a Heegaard surface and $M$ is not a $3$-sphere. 
So one of the circles in $\partial Q \cap \Sigma_y$ 
bounds a compressing disk for $\Sigma_y$. 
Define $D_y$ as such a disk. 
By definition, $D_y \cap D_{y_j}=\emptyset$ for $1 \le j \le 2^{i+1}-2$ and 
$\mathrm{diam}\,\partial D_y < 3^{1+\cdots+(i+1)} \cdot \delta$, 
which proves the claim. 
\end{proof}

\subsection*{Step~5: The conclusion} 
We now finish the proof of Theorem~\ref{thm:TIP_minimal_surface}. 
If $(0,0)$ and $(2,0)$ can be connected by a path in $[0,2] \times B$ 
without meeting $U$, 
it defines an isotopy $h'_t$ with the desired property, 
proving the theorem. 
Thus, it suffices to show that $[0,2] \times B \setminus U$ is path-connected. 
We will prove this by contradiction. 
Recall that $S=\partial B^d_{1-\epsilon}$. 

\begin{claim}
There exists a compact orientable $d$-manifold $Y$ 
in $U$ with $\partial Y=0 \times S$. 
\end{claim}

\begin{proof}
Consider the map $f: [0,2] \times \mathrm{int}\,B \to \mathbb{R}$ 
given by $f(t,s):=\gamma(\Sigma_{ts})$. 
Since $\Sigma_{ts}$ varies smoothly for $(t,s) \in [0,2] \times \mathrm{int}\,B$, 
$f$ is a continuous function. 
By the smooth approximation theorem, $f$ is approximated by a smooth map $f'$.   
Let $r \in \mathbb{R}$ be a regular value of $f'$ just below $\delta$. 
By assumption, one of the components of $f'^{-1}(r)$, say $Y'$, 
separates $(0,0)$ from $(2,0)$.  
On the other hand, by construction, 
if $(t,s)$ is close enough to $[0,2] \times \partial B$, 
then $f'(t,s) <r$.  
This implies that $\partial Y' \subset \{0,2\} \times B$. 
By Lemmas~\ref{lem:star} and \ref{lem:extension}~(3), 
$Y'$ extends to a $d$-manifold $Y$ in $U$ 
with $\partial Y=0 \times S$. 
\end{proof}

Pick a triangulation $\mathcal{Y}$ of $Y$ such that 
the diameter of any simplex of $\mathcal{Y}$ is small enough. 
Let $\mathcal{Y}'$ be the barycentric subdivision of $\mathcal{Y}$. 
By Lemma~\ref{lem:s-map}, 
we can find a simplicial map $\bar{\psi}:\mathcal{Y}' \to \Gamma(\Sigma)$. 
By Lemma~\ref{lem:h-nontrivial}, 
the restriction of $\bar{\psi}$ on $\partial Y$ must be homotopic to $\varphi$. 
Thus, $\varphi$ is homologically trivial. 
If $d \neq 2$, 
the Hurewicz theorem implies that $\varphi$ is homotopically trivial, 
contradicting the choice of $\varphi$. 

If $d=2$, we can deduce a contradiction as follows.  
Let $V$ and $W$ be the handlebodies in $M$ bounded by $\Sigma$. 
Denote by $\Gamma_V(\Sigma)$ (resp. $\Gamma_W(\Sigma)$) 
the subcomplex of $\Gamma(\Sigma)$ 
spanned by compressing disks for $\Sigma$ that lie in $V$ 
(resp. $W$). 
Furthermore, denote by $\Gamma_{VW}(\Sigma)$ 
the union of all simplices that contain vertices in both 
$\Gamma_V(\Sigma)$ and $\Gamma_W(\Sigma)$. 
Thus, $\Gamma(\Sigma)=
\Gamma_V(\Sigma) \cup \Gamma_{VW}(\Sigma) \cup \Gamma_W(\Sigma)$. 
Recall that the choice of $\varphi$ is arbitrary as long as it is homotopically nontrivial. 
By Claim~2.7 in \cite{Bac10}, we can assume that 
$\varphi$ is represented by a loop $\gamma$ in $\Gamma(\Sigma)$ with 
the following properties: 
\begin{enumerate}
\item[(a)] 
$\gamma$ can be expressed as $e \cup \gamma_V \cup e' \cup \gamma_W$, 
where $e$, $e'$ are edges in $\Gamma_{VW}(\Sigma)$ 
while $\gamma_V$, $\gamma_W$ are paths in 
$\Gamma_V(\Sigma)$ and $\Gamma_W(\Sigma)$, respectively. 
\item[(b)] 
$e$ is in the different component of $\Gamma_{VW}(\Sigma)$ from $e'$.  
\end{enumerate}

By definition, $\bar{\psi}(\partial Y)=\gamma$. 
Note that 
$\bar{\psi}^{-1}(\Gamma_V(\Sigma)) \cap \bar{\psi}^{-1}(\Gamma_W(\Sigma))
=\emptyset$. 
This along with (a) implies  
that there exists an arc in $\bar{\psi}^{-1}(\Gamma_{VW}(\Sigma))$ 
connecting $\bar{\psi}^{-1}(e)$ and $\bar{\psi}^{-1}(e')$. 
This contradicts (b) and  
completes the proof of Theorem~\ref{thm:TIP_minimal_surface}. 
\qed

\subsection*{Proof of Corollary~\ref{cor:MCG_minimal surface}}
We can now prove Corollary~\ref{cor:MCG_minimal surface}. 
By Theorems~\ref{thm:finitely_generatedness} and \ref{thm:TIP_minimal_surface}, 
$\pi_1(\mathcal{H}(M,\Sigma))$ is finitely generated, 
and this group projects onto $\mathrm{Isot}(M,\Sigma)$, 
the subgroup of $\MCG(M,\Sigma)$ 
that consists of maps $(M,\Sigma) \to (M,\Sigma)$ isotopic to $\mathrm{id}_M$. 
Since $\MCG(M)$ is a finite group, $\mathrm{Isot}(M,\Sigma)$ has finite index in 
$\MCG(M,\Sigma)$. 
Thus, $\MCG(M,\Sigma)$ is also finitely generated. 

\section{Examples: Heegaard splittings of $({\rm surface}) \times I$}
\label{sec:example}

Let $F$ be a closed orientable surface of genus $g \ge 2$, 
and set $P:=F \times [-1,1]$. 
For $n \in \mathbb{N}$, 
we consider the Heegaard surface $\Sigma^n=\Sigma$ for $P$ constructed as follows. 
Let $-1<r_1<\cdots <r_{n+1}<1$ and $F_i:=F \times \{r_i\}$. 
Fix a vertical arc $a_i$ in $F \times [r_i,r_{i+1}]$ connecting $F_i$ to $F_{i+1}$ and 
let $N(a_i) \cong B^2_1 \times [r_i,r_{i+1}]$ be 
a neighborhood of $a_i$ in $F \times [r_i,r_{i+1}]$. 
Then, $N(a_i)$ intersects $F_i$ (resp. $F_{i+1}$) in a disk $O^-_i:=B^2_1 \times \{r_i\}$ (resp. $O^+_i:=B^2_1 \times \{r_{i+1}\}$). 
Define $\Sigma$ as the surface obtained from $\bigcup F_i$ by 
replacing $\bigcup O^\pm_i$ with the tubes $\bigcup \partial B^2_1 \times [r_i,r_{i+1}]$. 
Then, $\Sigma$ cuts $P$ into two compact $3$-manifolds $V$ and $W$, 
each obtained from $F \times I$ and handlebodies 
by connecting them with $1$-handles.  
More precisely, if $n$ is even  
\[V \cong W = F \times [-1,r_1]  
\bigcup_{i~even} (F_i \setminus \mathrm{int}\,O^-_i) \times [r_i,r_{i+1}] 
\bigcup_{i~odd} N(a_i).\]
The odd case is the same except that $V$ is a handlebody 
and $\partial W \cap \partial P=\partial P$.  
Thus, $\Sigma$ is a Heegaard surface of genus $(n+1)g$. 
Furthermore, as a consequence of the main theorem of Lee \cite{Lee15}, we have 

\begin{theorem}[\cite{Lee15}]\label{thm:lee}
The topological index of $\Sigma^n$ is at most $n$. 
\end{theorem}

Note that when $n>1$, $\Sigma^n$ is stabilized 
and hence the index of $\Sigma^n$ is at least $2$. 
(The precise index is not known though it is likely equal to $n$.) 

In this section, 
we see that the proof in the previous section applies to 
this example after slight modification:  
We will show  

\begin{theorem}\label{thm:product}
For every $n$, $(P,\Sigma^n)$ satisfies the thick isotopy property. 
\end{theorem}

The setting of the previous section is different from this example 
in the following two points: 
1) $P$ has boundary; 
and 2) we do not know if $\pi_{d-1}(\Gamma(\Sigma))$ is finitely generated or not,
where $d>0$ denotes the index of $\Sigma$. 

1) is concerned with Theorem~\ref{thm:min-max}. 
Although the theorem assumes that the $3$-manifold is closed, 
the same conclusion also holds 
when the $3$-manifold has a boundary whose mean curvature field is inward pointing 
(cf. \cite[Theorem~2.1]{MN12}). 
Furthermore, we can endow $P$ with a hyperbolic metric so that
$\partial P$ satisfies this condition. 
Indeed, $F$ can be embedded in some hyperbolic $3$-manifold 
as a strictly stable minimal surface, 
and we can identify $P$ 
with a small tubular neighborhood of $F$. 
Thus, 1) is not the issue. 

As for 2), we need to show Lemma~\ref{lem:extension}
without the assumption on $\pi_{d-1}(\Gamma(\Sigma))$.   
The idea of the proof is as follows: 
Given a sweep-out $\{\Sigma_s\}_{s \in B}$,  
we can squeeze $\{\Sigma_s\}_{s \in B}$ 
into a thin product region by an isotopy $q_t$ 
so that every slice has bounded area. 
Using this isotopy we can construct an extension of $\{\Sigma_s\}_{s \in B}$ 
with appropriate properties. 

In what follows, we assume that $n$ is even for simplicity, although the same argument applies to an odd $n$ with straightforward modifications.  

\begin{figure}
\centering
\begin{overpic}[width=7cm,clip]{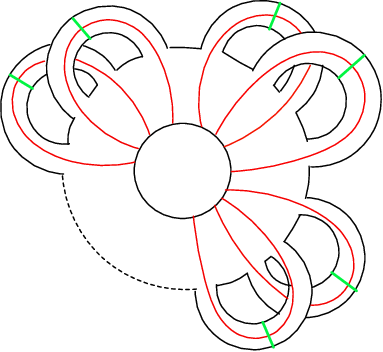}
\put(80,78){$\partial O^-_2$}
\put(30,40){$\partial O^+_1$} 
\end{overpic}
\caption{Arcs in $F_2 \setminus (O^+_1 \cup O^-_2)$: 
$c^1_1 \times \{r_2\},\ldots,c^1_{2g} \times \{r_2\}$ (green) and 
$e^2_1 \ldots,e^2_{2g}$ (red).}
\label{fig:arc_sys}
\end{figure}

We begin with a few definitions. 
For $1 \le i \le n+1$, put  $F_i^\circ:=F_i \setminus \Int O^-_i$. 
We say a compressing disk $D \subset F_i^\circ \times [r_i,r_{i+1}]$ is {\em vertical} 
if $D=c \times [r_i,r_{i+1}]$ for some properly embedded arc $c$ in $F_i^\circ$. 
Let $\{c^1_j\}^{2g}_{j=1}$ be properly embedded, pairwise disjoint arcs in $F_1^\circ$, 
which cut $F_1^\circ$ into a disk. 
Similarly, let $\{e^2_j\}^{2g}_{j=1}$ be arcs in $F_2^\circ$ satisfying the same condition. 
Furthermore, we may choose these arcs so that 
$|(c^1_j \times \{r_2\}) \cap e^2_{j'}|=\delta_{jj'}$, 
as depicted in Figure~\ref{fig:arc_sys}. 
Repeating the same argument for $3 \le i \le n+1$, 
we obtain collections of arcs $\{c^i_j\}^{2g}_{j=1}$ and $\{e^i_j\}^{2g}_{j=1}$. 
These arcs define vertical disks 
$C^i_j:=c^i_j \times [r_{i-1},r_i]$, $E^i_j:=e^i_j \times [r_i,r_{i+1}]$ in $V$, $W$ respectively. 
By definition, $|C^i_j \cap E^{i+1}_{j'}|=\delta_{jj'}$. 
Set $\mathscr{C}:=\{C^i_j\}$ and  
$\mathscr{E}:=\{E^i_j\}$. 

Suppose 
that $D$ and $D'$ are (possibly isotopic) disjoint compressing disks 
on the same side of $\Sigma$, say $V$. 
Let $c \subset \partial V$ be a simple arc intersecting $D \cup D'$ in its endpoints. 
We say $D''$ is obtained from $D$ and $D'$ by a {\em band summing} along $c$ 
if $D''$ is obtained from $D \cup D'$ 
by replacing $N(c) \cap (D \cup D')$ 
with a rectangle $b=\mathrm{fr}(N(c)) \setminus (D \cup D')$. 
Here $N(c)$ is a neighborhood of $c$ in $V$. 
We call such a rectangle $b$ a {\em band}.

\begin{lemma}\label{lem:band-sum}
Any compressing disk $D$ in $V$ (resp. $W$) is obtained 
from disks in $\mathscr{C}$ (resp. $\mathscr{E}$)
by a sequence of band summing.  
\end{lemma}

\begin{proof}
The proof is by induction 
on the intersection number between $D$ and $\mathscr{C}$. 
First, suppose that $D$ is disjoint from $\mathscr{C}$. 
Note that $\mathscr{C}$ cuts $V$ 
into $F \times I$ with ``scars'', 
each of which corresponds to a foot of a $1$-handle dual to a disk in $\mathscr{C}$. 
So, we can view $D$ as a disk in $F \times I$ and 
$D$ is isotopic to a disk $D'$ on $\partial (F \times I)$ 
that contains some scars.  
In other words, $D'$ is a neighborhood of 
the union of some scars and arcs connecting them. 
Thus, $D$ is the result of band summing of disks in 
$\mathscr{C}$. 

Next, suppose that $|D \cap \mathscr{C}|>0$. 
By the innermost disk argument, we may assume that $D$ intersects 
$\mathscr{C}$ only in finitely many arcs. 
Every arc in $D \cap \mathscr{C}$ cobounds a bigon 
$\Delta \subset \mathscr{C}$ 
together with a subarc in $\bigcup \partial \mathscr{C}$. 
Choose a bigon $\Delta$ such that $\mathrm{int}\, \Delta \cap D =\emptyset$.  
(In other words, $\Delta$ is an outermost bigon.)
$\partial$-compressing $D$ along $\Delta$ yields the new compressing disks $D'$, $D''$ 
which, by induction, are obtained from disks in $\mathscr{C}$
by a sequence of band summing.  
As $D$ is obtained by band summing from $D'$ and $D''$, 
the conclusion follows. 
\end{proof}

For $t \in I$, 
let $q_t:F \times [-1,1] \to F \times[-1,1]$ be the map given by $q_t(x,r)=(x,(1-t)r)$. 

\begin{lemma}\label{lem:disk-squeezing} 
Suppose that 
$\{D_v\}$ is a finite collection of compressing disks for $\Sigma$.   
For every $\epsilon>0$, 
there exists $\epsilon'>0$ satisfying the following. 
For every $v$, there exists a disk $D'_v$ isotopic to $D_v$ 
such that for any $t \in (1-\epsilon',1)$
\[\mathrm{Area}(q_t(D'_v)) < \epsilon.\] 
Moreover, if $D_v$ and $D_w$ are disjoint,  
then the same is true for $D'_v$ and $D'_w$. 
\end{lemma}

\begin{proof} 
Note that if $D_v$ is a vertical disk, then $\mathrm{Area}(q_t(D_v)) \to 0$ 
as $t \to 1$.
The idea of the proof is rather simple: 
By Lemma~\ref{lem:band-sum},  
$D_v$ can be expressed as a band sum of vertical disks,  
and thinning each band isotopes $D_v$ to a disk $D'_v$ 
such that $\mathrm{Area}(q_t(D'_v)) \to 0$. 
We modify this argument and 
show that  $D'_v$,$D'_w$ can be taken to be disjoint if $D_v$,$D_w$ are disjoint.   

We start with some setups. 
Let $\#C^i_j$ be a disk in $V$ obtained from $C^1_j,C^3_j,\ldots,C^i_j$  
by band summing along arcs in $\partial E^2_j,\partial E^4_j,\ldots, \partial E^{i-1}_j$, 
as shown in Figure~\ref{fig:disk_sys}. 
\begin{figure}
\centering
\begin{overpic}[height=3cm,width=7cm,clip]{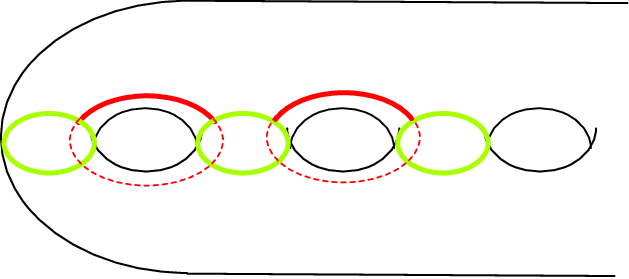}
\end{overpic}
\caption{$\#C^i_j$ is obtained from $C^1_j,C^3_j,\ldots,C^i_j$ (green) 
by band summing along arcs in $\partial E^2_j,\partial E^4_j,\ldots, \partial E^{i-1}_j$ 
(bold red).}
\label{fig:disk_sys}
\end{figure}
Set $\mathscr{C}^{\#}:=\{\#C^i_j\}$. 
Note that $\mathscr{C}^{\#}$ is orthogonal to $\mathscr{E}$ 
(i.e. $|\#C^i_j \cap E^{i+1}_j|=1$ and $|\#C^i_j \cap E^{i'}_{j'}|=0$ otherwise), 
and from this point of view $\mathscr{C}^{\#}$ is easier to handle than $\mathscr{C}$.  
So, we work with $\mathscr{C}^{\#}$ in the following argument. 

Observe that $C^i_j$ is recovered from $\#C^i_j$ and $\#C^{i-2}_j$ 
by taking band summing. 
By Lemma~\ref{lem:band-sum}, 
$D_v$ can be written in this form: 
\[D_v= \bigcup^{m_v}_{p=1} F^p_v \cup \bigcup^{n_v}_{k=1} b^k_v,\]
where $F^p_v$ is a $\mathscr{C}^{\#}$- or $\mathscr{E}$-component of $D_v$, 
depending on if $D_v \subset V$ or $D_v \subset W$ 
(that is, $F^p_v$ is isotopic to a disk in $\mathscr{C}^{\#}$ or $\mathscr{E}$ 
after attaching bigons along arcs $F^p_v \cap \bigcup^{n_v}_{k=1} b^k_v$),  
and $b^k_v$ is a band. 
By taking the bands connecting $C^i_j$'s to be  thin enough, we may assume that 
\begin{equation}\label{eq:thin-band}
\Area(q_t(\#C^i_j)) < \frac{\epsilon}{4m_v} 
\end{equation}
for all $v$ and $t$ sufficiently close to $1$. 
Similarly, 
\begin{equation}\label{eq:vertical_disk}
\Area(q_t(E^i_j)) < \frac{\epsilon}{4m_v} 
\end{equation}
for all $v$ and $t$ sufficiently close to $1$. 

We define the disk $D'_v$ isotopic to $D_v$ as follows.   
Let $N(\mathscr{C}^{\#}) \cong \mathscr{C}^{\#} \times [-1,1]$, 
$N(\mathscr{E}) \cong \mathscr{E}\times [-1,1]$ 
be small product neighborhoods of $\mathscr{C}^{\#}$ and $\mathscr{E}$, respectively. 
Isotope $D_v$ so that $F^p_v$ is a subdisk of 
$\mathscr{C}^{\#} \times \{u\}$ or $\mathscr{E} \times \{u\}$ for some $u$,
depending on if $D_v \subset V$ or $D_v \subset W$. 
Let $F'^p_v$ denote the resulting disk that corresponds to $F^p_v$. 
If $b^k_v$ is a band of $D_v$, 
there is a unique rectangle $\ol{b}^k_v \subset \Sigma$ 
(the ``shadow'' of $b^k_v$)  
determined by the band sum structure of $D_v$.  
\begin{claim*}
There exists a point  $x \in \Sigma \setminus (N(\mathscr{C}^{\#}) \cup N(\mathscr{E}))$ 
such that for every $v$ and $k$, $\ol{b}^k_v$ does not contain $x$.  
\end{claim*}
\begin{proof}
We can see this, for example, as follows. 
Let 
$S \subset \Sigma \setminus (N(\mathscr{C}^{\#}) \cup N(\mathscr{E}))$ 
be a three-holed sphere bounded by essential simple closed curves. 
Suppose that $D_v$ intersects $S$ minimally. 
Then, every rectangle $\ol{b}^k_v$ intersects $S$ 
in finitely many pairwise disjoint rectangles. 
Note that there are only six possible types for such rectangles 
up to twisting around $\partial S$, 
as shown in Figure~\ref{fig:bd-twist} (a). 
Consider components $\ol{b}_1$, $\ol{b}_2$ of $\ol{b}^k_v \cap S$, $\ol{b}^l_w \cap S$ respectively. 
If $\partial S$, $\ol{b}_1$ and $\ol{b}_2$ form a triangle, 
we can push it off $S$ as in Figure~\ref{fig:bd-twist} (b). 
By repeating this operation until all such triangles are eliminated, 
we may assume that $\ol{b}_1$ and $\ol{b}_2$ do not intersect near $\partial S$ 
if $\ol{b}_1$, $\ol{b}_2$ are of different types. 
This implies $\partial S \setminus \bigcup_{k,v} \ol{b}^k_v \neq \emptyset$. 
\begin{figure}
\centering
\begin{overpic}[width=10cm,height=2.5cm,clip]{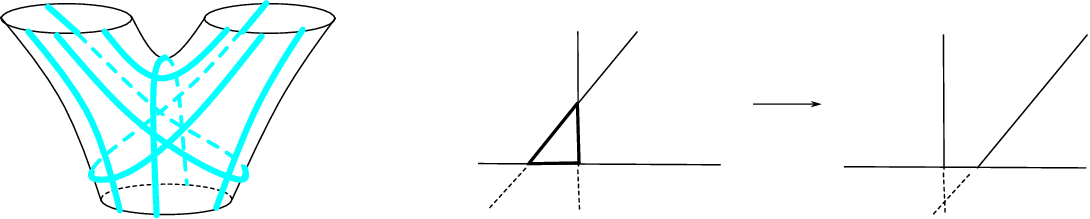}
\put(35,-20){(a)}
\put(200,-20){(b)} 
\put(110,15){$\partial S$}
\put(140,60){$\ol{b}_1$}
\put(170,60){$\ol{b}_2$}
\end{overpic}
\vspace{2em}
\caption{(a) There are only six possible types for rectangles in $S$ up to twisting around $\partial S$. 
	   (b) A triangle formed by $\partial S$ and sides of $\ol{b}_1$, $\ol{b}_2$ 
	   can be pushed off $S$.}
\label{fig:bd-twist}
\end{figure}
\end{proof}

Let $O_x$ be a small neighborhood $O_x \cong \Int B^2_1$ of $x$ in $\Sigma$. 
Imagine that 
we expand $O_x$ across $\Sigma$ and 
push $\ol{b}^k_v$ simultaneously into a thin neighborhood of a graph $G \subset \Sigma$ 
depicted in Figure~\ref{fig:thickened-graph}: 
The graph $G$ can be taken so that it satisfies the following: 
\begin{itemize}
\item $G$ is a deformation retract of $\Sigma \setminus O_x$.   
\item $G$ contains the boundaries of disks 
	$\#C^1_1, E^2_1, \#C^3_1, E^4_1, \ldots, \#C^{n-1}_{2g},E^n_{2g}$  
	as its loops. 
\end{itemize}
\begin{figure}
\centering
\begin{overpic}[width=7cm,clip]{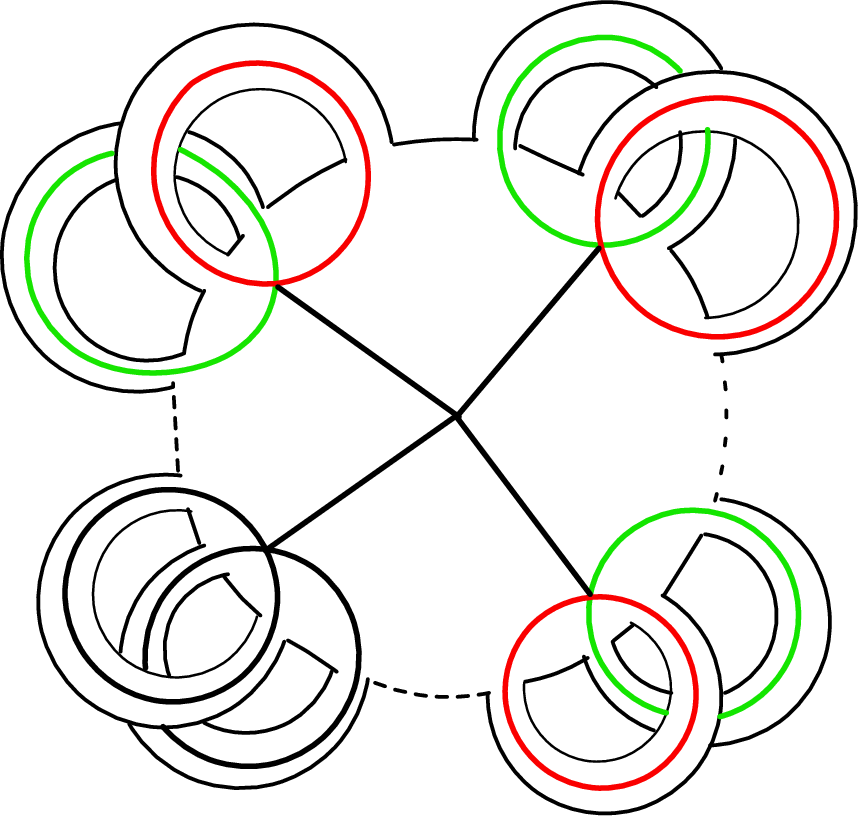}
\put(100,75){$G$}
\put(50,100){$\#C_1^1$} 
\put(80,120){$E^2_1$}
\put(106,125){$\#C^3_1$}
\put(136,105){$E^4_1$}
\put(135,78.3){$\#C^{n-1}_{2g}$}
\put(100,40){$E^n_{2g}$}
\end{overpic}
\caption{$G$ is a deformation retract of $\Sigma \setminus O_x$.}
\label{fig:thickened-graph}
\end{figure}
This map extends to an ambient isotopy 
that leaves $N(\mathscr{C}^{\#})$ and $N(\mathscr{E})$ invariant. 
Isotope further each $b^k_v$ to a rectangle $b'^k_v$ near $\ol{b}_v^k$ so that 
\begin{equation}\label{eq:thin-band2}
\sum^{n_v}_{k=1} \Area(q_t(b'^k_v)) <\frac{\epsilon}{2}
\end{equation} 
for $t$ sufficiently close to $1$. 
Let $D'_v$ be the resulting disk. 

We see that $D'_v$ satisfies the desired property.  
If $D_v \subset V$, 
by taking $F'^p_v$ to be close enough to $\#C^i_j$, we may assume that 
\[\Area(q_t(F'^p_v)) < \Area(q_t(\#C^i_j)) + \frac{\epsilon}{4m_v}\] 
Similarly,   
\[\Area(q_t(F'^p_v)) < \Area(q_t(E^i_j)) + \frac{\epsilon}{4m_v}\] 
if $D_v \subset W$. 
By (\ref{eq:thin-band}) or (\ref{eq:vertical_disk}) we have 
\[\Area(q_t(F'^p_v)) < \frac{\epsilon}{2m_v}\] 
for $t$ sufficiently close to $1$. 
Combining this with (\ref{eq:thin-band2}), we have 
\[\Area(q_t(D'_v)) < \sum^{m_v}_{p=1} \Area(q_t(F'^p_v))+\sum^{n_v}_{k=1} \Area(q_t(b'^k_v)) < \epsilon\]  
for $t$ sufficiently close to $1$. 

Finally, we see that if $D_v$, $D_w$ are disjoint, 
then so are $D'_v$, $D'_w$. 
Observe that if $F^p_v$ is already in an appropriate position,  
then we can isotope $D_w$ in the complement of $\bigcup_p F^p_v$ 
so that the same is true for $F^q_w$. 
The isotopies from $b^k_v$ to $b'^k_v$ and from $b^l_w$ to $b'^l_w$ 
can be done simultaneously. 
Thus, we can take $F'^p_v$,$F'^q_w$,$b'^k_v$ and $b'^l_w$ so that  
$F'^p_v \cap F'^q_w=b'^k_v \cap b'^l_w=\emptyset$, which proves the lemma. 
\end{proof}

Let $\delta >0$ as in Section~\ref{sec:proof_TIP_minimal_surface}. 
Set 
\[\delta':=\min \left\{\delta, \inf_{t \in [0,1)} \gamma(q_t(\Sigma))\right\},\]
and 
\[C':=\max_{t \in I} \mathrm{Area}(q_t(\Sigma))+1.\] 
Note that $\delta'>0$ by Lemma~\ref{lem:band-sum}. 

\begin{lemma}\label{lem:squeezing}
Suppose that $\psi:S^{d-1} \to \Gamma(\Sigma)$ is given,  
and $\{\Sigma_s\}_{s \in B}$ is the sweep-out constructed from $\psi$ 
as in Step~1 of Section~\ref{sec:proof_TIP_minimal_surface}. 
Let $t' \in (0,1)$. 
Then, there exists an extension $\{\Sigma_{ts}\}_{(t,s) \in [0,t'] \times B}$ 
of $\{\Sigma_s\}_{s \in B}$ such that 
\begin{enumerate}
\item $\Sigma_{t0}=\Sigma$ for $t \in [0,t']$.  
\item $\mathrm{Area}(\Sigma_{t's}) <C'$ for $s \in B$.  
\item $B \setminus U_{t'}$ is a star-shaped region, 
where $U_{t'}$ is the set of those points $s$ 
such that $\Sigma_{t's}$ is $\delta'$-compressible. 
\end{enumerate}
\end{lemma}

\begin{proof}
Fix $t \in I$ for a moment. 
Let $\ell_t:B_{1/2t} \to [0,t]$ be the map 
given  by $\ell_t(s)=2|s|$.  
Set $B_{[1/2t,1]}=B \setminus \mathrm{int}\,B_{1/2t}$.
Let $\varpi_t: B_{[\frac{1}{2}t,1]} \to B$ be the map given by 
\[\varpi_t(rs)=\frac{2r-t}{2-t}s\] 
for $r \in [1/2t,1]$ and $s \in S^{d-1}=\partial B$. 
Then, define the sweep-out $\{\Sigma_{ts}\}_{s \in B}$ by 
\[\Sigma_{ts}:= 
\begin{cases}
q_{\ell_t(s)}(\Sigma) & s \in B_{\frac{1}{2}t},\\
q_t\left(\Sigma_{\varpi_t(s)}\right) & s \in B_{[\frac{1}{2}t,1]}.
\end{cases}
\]
As $\Sigma_{ts}$ varies smoothly with $t$, 
we obtain the sweep-out $\{\Sigma_{ts}\}_{(t,s) \in [0,1) \times B}$. 
It follows from the construction that $\Sigma_{t0}=\Sigma$ for $t \in I$. 

It follows from the definition of $\delta'$ 
that $B_{\frac{1}{2}t} \cap U_t=\emptyset$ for all $t$. 
On the other hand, after relabeling via $\varpi_t: B_{[\frac{1}{2}t,1]} \to B$, 
$\{\Sigma_{ts}\}_{s \in B_{[1/2t,1]}}$ 
is nothing but the image of $\{\Sigma_s\}_{s \in B}$ by $q_t$. 
Thus, the item (3) follows from Lemma~\ref{lem:star}. 

It remains to see (2) holds. 
We see that for $t$ sufficiently close to $1$, 
$\{\Sigma_{ts}\}_{s \in B}$ can be modified 
so that $\{\Sigma_{ts}\}_{s \in B}$ satisfies (2). 

For each vertex $v$ of $S^{d-1}$, 
let $D_v$ be a representative of $\psi(v)$ 
such that if $v$ and $w$ are in the same simplex, 
$D_v$ and $D_w$ are disjoint from each other. 
By Lemma~\ref{lem:disk-squeezing}, 
there exist $\epsilon'>0$ and a disk $D'_v$ isotopic to $D_v$ 
such that 
\begin{equation}\label{eq:disk}
\mathrm{Area}(q_t(D'_v)) < \frac{1}{4N}
\end{equation}
for $t \in (1-\epsilon',1)$. 
Here, $N>0$ is the number of vertices of $S^{d-1}$.  
Fix $t  \in (1- \epsilon',1)$.  

Consider the sweep-out $\Sigma'_{s \in B}$ given as follows. 
For each $v$, 
set $D''_v:=q_t\left(D'_v\right)$ and 
fix an identification $\iota_v:B^2_1 \times [-1,1] \to N(D''_v)$. 
By Lemma~\ref{lem:disk-squeezing}, 
$D''_v$ and $D''_w$ are disjoint from each other 
if $v$, $w$ are contained in the same simplex of $S^{d-1}$. 
So, we can construct a sweep-out $\{\Sigma''_s\}_{s \in B}$ with $\Sigma''_0=q_t(\Sigma)$ 
from $D''_v$, $N(D''_v)$ and $\iota_v$ 
as in Step~1 of Section~\ref{sec:proof_TIP_minimal_surface}. 
As before, define 
\[\Sigma'_s:= 
\begin{cases}
q_{\ell_t(s)}(\Sigma) & s \in B_{\frac{1}{2}t},\\
q_t\left(\Sigma''_{\varpi_t(s)}\right) & s \in B_{[\frac{1}{2}t,1]}.
\end{cases}
\] 

By definition, 
$\mathrm{Area}(\Sigma'_s)<C'$ for $s \in B_{1/2t}$. 
Thus, it suffices to see that $\mathrm{Area}(\Sigma'_s)<C'$ for $s \in B_{[1/2t,1]}$. 

On each simplex $\sigma$ of $B$, we can express the area of $\Sigma'_s$ as  
\[\mathrm{Area}(\Sigma'_s) = 
\mathrm{Area}\left(\Sigma'_s \setminus \bigcup_{v \in \sigma} N(D''_v)\right)
+\sum_{v \in \sigma} \mathrm{Area}\left(\Sigma'_s \cap N(D''_v)\right).
\]

For $\tau \in I$, 
let $\varsigma_\tau:B^2_1 \times [-1,1] \to B^2_1 \times [-1,1]$ be the map given by 
$\varsigma_\tau(x,u):=(x,(1-\tau)u)$. 
Set $\lambda^v_\tau:=\iota_v \circ \varsigma_\tau \circ \iota^{-1}_v$. 
Then, $\lambda^v_\tau$ shrinks $N(D''_v)$ in the $I$-direction.  
Letting $\tau \to 1$,  we have  
\[\mathrm{Area}(\lambda^v_\tau(\Sigma'_s \cap N(D''_v))) 
< 2\mathrm{Area}(D''_v)+\frac{1}{2N},\]
for $s \in \sigma \cap B_{[1/2t,1]}$ and $\tau$ sufficiently close to $1$. 
Thus, after shrinking each $N(D''_v)$ by $\lambda^v_\tau$, we have 
\begin{align*}\label{eq:shrinking}
\mathrm{Area}(\Sigma'_s) 
&< \mathrm{Area}\left(\Sigma'_s \setminus \bigcup_{v \in \sigma} N(D''_v)\right) 
+\sum_{v \in \sigma} 2\mathrm{Area}(D''_v)+\sum_{v \in \sigma} \frac{1}{2N}\\
&<\mathrm{Area}\left(\Sigma'_s \setminus \bigcup_{v \in \sigma} N(D''_v)\right)+1 
\end{align*}
on $\sigma \cap B_{[1/2t,1]}$. 
Here, we used the inequality \ref{eq:disk} for the second line.  
By construction, $\Sigma'_s$ coincides with $q_t(\Sigma)$ 
in the complement of $\bigcup N(D''_v)$,   
and hence 
$\mathrm{Area}(\Sigma'_s \setminus \bigcup N(D''_v))
=\mathrm{Area}(q_t(\Sigma) \setminus \bigcup N(D''_v))$. 
Combining this with the above inequality implies 
\[\mathrm{Area}(\Sigma'_s)< C'\] 
for $s \in B_{[1/2t,1]}$. 

By the same argument at the end of the proof of Lemma~\ref{lem:extension}, 
we can interpolate between 
$\{\Sigma_{ts}\}_{s \in B}$ and $\{\Sigma'_s\}_{s \in B}$ via isotopy. 
The resulting sweep-out satisfies (2) 
as well as the other two conditions. 
This completes the proof. 
\end{proof}

Having established Lemma~\ref{lem:squeezing}, 
Theorem~\ref{thm:product} now follows 
by the same argument as the proof of Theorem~\ref{thm:TIP_minimal_surface}: 
we have only to use Lemma~\ref{lem:squeezing} instead of Lemma~\ref{lem:extension}. 
Also, a similar argument to Corollary~\ref{cor:MCG_minimal surface} proves the following.  

\begin{theorem}\label{thm:MCG_surf_x_I}
For every $n \in \mathbb{N}$, 
$\mathrm{MCG}(P, \Sigma^n)$ is finitely generated. 
\end{theorem}

\appendix\section{The disk complex of a Heegaard splitting of $S^3$}\label{sec:appA}

In this appendix, 
we prove that the disk complex of a genus $\ge 2$ 
Heegaard surface $\Sigma$ of $S^3$ is not of finite type, 
as mentioned in Section~\ref{sec:intro}.  
More precisely, we show:

\begin{proposition}\label{prop:infinite}
Both $H_{2g-2}(\mathcal{C}(\Sigma))$ and 
$H_{2g-2}(\Gamma(\Sigma))$ 
are not finitely generated. 
\end{proposition}

Here, $\mathcal{C}(\Sigma)$ is the {\em curve complex} of $\Sigma$ 
defined as follows. 
The vertices of $\mathcal{C}(\Sigma)$ 
are isotopy classes of essential simple closed curves on $\Sigma$. 
A $k$-tuple of vertices spans a $k$-simplex if these vertices admit 
pairwise disjoint representatives.  
We first prove the proposition for the curve complex, 
and then for the disk complex.  
(The assertion for the curve complex is previously known as Theorem~1.4 of Ivanov-Ji \cite{IJ}. 
But here we reprove this fact in a way that can apply to the disk complex.)  

For brevity, set $H:=H_{2g-2}(\mathcal{C}(\Sigma))$. 
We give the norm for $H$ as follows:   
For any simplicial chain $c=\sum a_i\sigma_i$, define 
\[\|c\|:=\sum |a_i|.\] 
For $\alpha \in H$, define 
\[\|\alpha\|:=\min_{[c]=\alpha} \|c\|.\] 
The next lemma follows from the definition. 

\begin{lemma}\label{lem:norm_invariant}
The norm $\| \cdot \|$ is invariant under the action of $\mathrm{MCG}(\Sigma)$.  
\end{lemma}

Take the tensor product $\mathbb{R} \otimes \widetilde{H}$. 
We think of $H$ as a subset of 
$\mathbb{R} \otimes H$. 

\begin{lemma}
There exists a norm on $\mathbb{R} \otimes H$ 
such that 
\[\|r \otimes \alpha\|=|r|\|\alpha\|\] 
for $r \in \mathbb{R}$ and $\alpha \in H$. 
\end{lemma}

\begin{proof}
We can define, for example, the norm on $\mathbb{R} \otimes H$
as follows. 
(The definition is similar to the {\em injective cross norm} first introduced in \cite{Gro53}.)   
Denote by $B_{H^*}$ the set of homomorphisms 
$\psi:H \to \mathbb{R}$ such that 
$|\psi(\alpha)| \le \|\alpha\|$ for all $\alpha \in H$. 
For $\beta=\sum r_i \otimes \alpha_i \in \mathbb{R} \otimes H$, 
define 
\[\|\beta\|:=\sup \left\{ \left| \sum r_i \psi(\alpha_i) \right|
\mid \psi \in B_{H^*} \right\}.\]

It follows from the definition that this is a seminorm and satisfies 
$\|r \otimes \alpha\|=|r|\|\alpha\|$  
for $r \in \mathbb{R}$ and $\alpha \in H$. 
If $\beta=\sum r_i \otimes \alpha_i$ and $\|\beta\|=0$, then 
\[|r_i|\|\alpha_i\| = 
\left|r_i\psi_{\alpha_i}(\alpha_i)\right| 
\le \left\|\sum r_i \otimes \alpha_i \right\|=0\] 
for all $i$. 
Here, $\psi_{\alpha_i} \in B_{H^*}$ is the homomorphism 
given by $\psi_{\alpha_i}(n\alpha_i)=n\|\alpha_i\|$ for $n \in \mathbb{Z}$ 
and $\psi_{\alpha_i}=0$ otherwise. 
Thus, $\beta=0$ and we conclude that this seminorm is in fact a norm. 
\end{proof} 

Now we prove 
that $H$ is not finitely generated. 
Denote by $\mathcal{B}(r_0) \subset \mathbb{R} \otimes H$ 
the ball of diameter at most $r_0>0$. 
Note that $H \neq 0$ by \cite[Theorem~1.3]{IJ} (or by \cite[Theorem~1.1]{Bro12}). 
By taking $r_0$ to be large enough, 
we may assume 
that $H \cap \mathcal{B}(r_0)$ 
contains a nontrivial homology class. 
Suppose, contrary to our claim, that 
$H$ is finitely generated.  
Then, $\mathbb{R} \otimes H$
is also finitely generated, 
and hence $\mathbb{R} \otimes H$ is equivalent to 
$\mathbb{R}^n$ with the Euclidean norm. 
It follows that 
$H \cap \mathcal{B}(r_0)$ is a finite set 
since any closed bounded set in 
$\mathbb{R} \otimes H \cong \mathbb{R}^n$ 
is compact.

It follows from Lemma~\ref{lem:norm_invariant} 
that $\MCG(\Sigma)$ maps $H \cap \mathcal{B}(r_0)$ to itself. 
Thus, for every point $\alpha \in H \cap \mathcal{B}(r_0)$, 
the stabilizer subgroup of $\alpha$ is an infinite group. 
But Corollary~5.3 in Irmer \cite{Irm23} says 
that every stabilizer subgroup for the action of
$\MCG(\Sigma)$ on $H \setminus \{0\}$ 
must be trivial or $\mathbb{Z}/2\mathbb{Z}$ 
(the latter case occurs only when $g=2$), 
a contradiction.  

Next, we see that $H_{2g-2}(\Gamma(\Sigma))$ is not finitely generated.
The proof is similar to the curve complex case. 
Note that the disk complex $\Gamma(\Sigma)$ can be identified 
with a subcomplex of $\mathcal{C}(\Sigma)$ in the following sense. 
Consider the map given by 
sending $[D] \in \Gamma(\Sigma)$ to $[\partial D] \in \mathcal{C}(\Sigma)$.  
(This is not an injection because $\partial D$ may bound a disk 
opposite to $D$.) 
This map is a homotopy equivalence 
between $\Gamma(\Sigma)$ and its image $\mathcal{D}_\Sigma$. 
Let $H_\mathcal{D}$ be the subspace of $H$ 
consisting of those homology classes $\alpha$ such that 
$\alpha$ is represented by a cycle in $\mathcal{D}_\Sigma$.  

In \cite{Bro12}, Broaddus identified a generator of 
$H$ as a $\mathbb{Z}\MCG(\Sigma)$-module, 
which is represented by a nontrivial $(2g-2)$-sphere in $\mathcal{C}(\Sigma)$. 
Surprisingly, such a sphere can be found within $\mathcal{D}_\Sigma$. 
See e.g. Figures~1-3 in \cite{CT20}. 
Thus, $H_\mathcal{D} \neq 0$. 
In particular, $H_\mathcal{D} \cap \mathcal{B}(r_0)$ contains a nontrivial homology class 
for $r_0$ sufficiently large. 

As before, suppose 
that $H_{2g-2}(\Gamma(\Sigma))$ is finitely generated. 
Then, so is $H_\mathcal{D}$ since $H_\mathcal{D}$ can be thought of as a subspace of 
$H_{2g-2}(\mathcal{D}_\Sigma) \cong H_{2g-2}(\Gamma(\Sigma))$.  
By the same argument as above, 
it follows that $H_\mathcal{D} \cap \mathcal{B}(r_0)$ is a finite set. 
Note that restricting maps $(S^3,\Sigma) \to (S^3,\Sigma)$ on $\Sigma$ defines 
the injection $\MCG(S^3,\Sigma) \to \MCG(\Sigma)$. 
So we can think of $\MCG(S^3,\Sigma)$ as a subgroup of $\MCG(\Sigma)$. 
By Lemma~\ref{lem:norm_invariant}, 
$\MCG(S^3,\Sigma)$ maps $H_\mathcal{D} \cap \mathcal{B}(r_0)$ to itself. 
As $\MCG(S^3,\Sigma)$ is an infinite group,  
the stabilizer subgroup of a point in $H_\mathcal{D} \cap \mathcal{B}(r_0)$ 
is an infinite group, 
again contradicting Corollary~5.3 in \cite{Irm23}. 
This completes the proof of Proposition~\ref{prop:infinite}. 
\qed 

As a consequence of Proposition~\ref{prop:infinite} and 
the main theorem of Appel \cite{App10} 
or Campisi-Torres \cite{CT20}, we have

\begin{theorem}\label{thm:infinite_bouquet}
Suppose that $\Sigma$ is a Heegaard surface for $S^3$ of genus $g \ge 2$. 
Then the disk complex $\Gamma(\Sigma)$ is 
homotopy equivalent to the bouquet of countably infinitely many spheres 
of dimension $2g-2$. 
\end{theorem}

\bibliographystyle{amsrefs}

\end{document}